\documentclass[10pt]{article}

\usepackage[english]{babel}
\usepackage[T1]{fontenc}
\usepackage[utf8]{inputenc}

\usepackage{amsfonts}
\usepackage{amsmath}
\usepackage{amssymb}
\usepackage[a4paper,left=3cm,right=3cm,bottom=4cm,top=3cm]{geometry}

\author{Vincent Jugé\footnote{LSV, CNRS \& ENS Cachan, Univ. Paris-Saclay, France}}
\title{Abelian Ramsey Length and Asymptotic Lower Bounds\footnote{This work is supported by ERC EQualIS (308087).}}
\date{}

\newcommand{\calA}{{\mathcal{A}}}

 \newcommand{\calV}{{\mathcal{V}}}

\newcommand{\bK}{{\textbf{K}}} 
 
 \newcommand{\bP}{{\textbf{P}}}
 
\newcommand{\bS}{{\textbf{S}}} \newcommand{\bT}{{\textbf{T}}}
\newcommand{\bU}{{\textbf{U}}} \newcommand{\bV}{{\textbf{V}}}
 
 \newcommand{\bZ}{{\textbf{Z}}}

\newcommand{\ba}{{\textbf{a}}}

 \newcommand{\bp}{{\textbf{p}}}

 \newcommand{\NN}{{\mathbb{N}}}
 
 \newcommand{\RR}{{\mathbb{R}}}

\newcommand{\ds}{\displaystyle}
\newcommand{\ol}{\overline}
\newcommand{\formula}[1]{\begin{center}$\ds #1$.\end{center}}
\newcommand{\formulaa}[1]{\begin{center}$\ds #1$,\end{center}}

\begin{document}
\maketitle

\medskip

This technical note aims at evaluating an asymptotic lower bound on \emph{abelian Ramsey lengths} obtained by Tao in~\cite{2014arXiv1406.0450T}.
We first provide the minimal amount of background necessary to define abelian Ramsey lengths,
and indicate the lower bound of Tao. We then focus on evaluating this lower bound.

\section{Introduction}

Let $\calA$ and $\calV$ be two alphabets. A word on $\calA$ is a finite sequence $\ba = a_1 a_2 \cdots a_k$ of elements of $\calA$.
The elements $a_i$ are called the letters of the word $\ba$, and the integer $k$ is the length of $\ba$.
For all elements $\alpha \in \calA$, we denote by $|\ba|_\alpha$ the cardinality of the set 
$\{i : a_i = \alpha\}$, i.e. number of occurrences of the letter $\alpha$ in the word $\ba$.
We also denote by $\calA^\ast$ the set of all words on $\calA$.
The words $a_i a_{i+1} \cdots a_j$ with $1 \leq i \leq j \leq k$, as well as the empty word, are called factors of $\ba$.

Consider now a word $\ba = a_1 a_2 \cdot a_k$ in $\calA^\ast$ and a word $\bp = p_1 p_2 \ldots p_\ell$ in $\calV^\ast$.
We say that $\ba$ contains $\bp$ in the abelian sense if
there exist non-empty words $\pi_1, \pi_2, \ldots, \pi_\ell$ in $\calA^\ast$ such that
the concatenated word $\pi_1 \pi_2 \ldots \pi_\ell$ is a factor of $\ba$, and such that,
for all integers $i, j$ and all letters $\alpha \in \calA$, if $p_i = p_j$, then $|\pi_i|_\alpha = |\pi_j|_\alpha$.
For instance, the word \emph{programmable} contains the word \emph{aab} in the abelian sense,
as can be seen by considering the words $\pi_1 =$ \emph{am}, $\pi_2 =$ \emph{ma} and $\pi_3 = $ \emph{ble}.

From this point on, we consider the infinite alphabet $\calV = \{v_i : i \in \NN\}$, where $\NN$ is the set of positive integers,
and we define the Zimin patterns $Z_i$ inductively by $Z_1 = v_1$ and $Z_{i+1} = Z_i v_{i+1} Z_i$.
It turns out that, for all integers $i, m \geq 1$ and all alphabets $\calA$ of cardinality $m$,
there exists an integer $L_{ab}(m,Z_i)$ such that
all words $\ba \in \calA^\ast$ with length at least $L_{ab}(m,Z_i)$ contain the word $Z_i$ in the abelian sense.

For all integers $m \geq 4$, Tao proves in~\cite{2014arXiv1406.0450T} that $L_{ab}(m,Z_i) \geq (1+\varepsilon_m(i)) \sqrt{\bK(m,i)}$ for all $i \geq 1$,
where $\varepsilon_m$ is a function such that $\lim_{+\infty}\varepsilon_m = 0$ and $\bK(m,i)$ is defined as
\begin{center}$\bK(m,i) = 2 \prod_{j=1}^{i-1} \bS(m,2^j)^{-1}$,\end{center}
where $\ds\bS(m,k) = \sum_{\ell = 1}^\infty \bT(m,k,\ell)$ and $\ds\bT(m,k,\ell) = \frac{1}{m^{k \ell}} \sum_{i_1+\ldots+i_m=\ell} \binom{\ell}{i_1~\ldots~i_m}^{k}$.

Yet, in order to obtain actual lower bounds on $L_{ab}(m,Z_i)$, it remains to evaluate the asymptotical behavior of $\bK(m,i)$.
We evaluate $\bK(m,i)$ up to a multiplicative constant that does not depend on $m$ or $i$.
More precisely, we prove the following inequalities, which hold for all $m \geq 4$ and $i \geq 1$:
\formula{ 2 \frac{m^{2^i}}{m^{i+1}} \geq \bK(m,i) \geq \frac{1}{21} \frac{m^{2^i}}{m^{i+1}} }

\pagebreak

\section{Auxiliary inequalities}

Before evaluating the lower bound $\bK(m,i)$, we prove a series of six inequalities that we will use subsequently.
We first study the function $f_x : y \mapsto y \ln(1+x/y)$ for $x,y > 0$.
An asymptotic evaluation proves that $\lim_{+\infty} f_x = x$.
Furthermore, we compute that $f''_x(y) = -\frac{x^2}{y(x+y)^2} < 0$ for $y > 0$.
It follows that $x > f_x(y)$ or, equivalenlty, that 
\begin{eqnarray}
(1+x/y)^y < e^x && \text{for all } x,y > 0.
\end{eqnarray}

We perform a similar study with the function $g : y \mapsto (y+1/2)\ln(1+1/y)$ for $y > 0$.
We find that $\lim_{+\infty} g = 1$ and that $g''(y) = \frac{1}{2y^2(y+1)^2} > 0$ for $y > 0$.
It follows that $g(y) > 1$ or, equivalently, that
\begin{eqnarray}
(1+1/y)^{y+1/2} > e && \text{for all } y > 0.
\end{eqnarray}

Again, we consider the function $h : y \mapsto 3 \ln(y) + \ln(2) - (y-1) \ln(2 \pi)$ for $y > 0$,
as well as the real constant $\lambda = \frac{4}{\pi^{3/2}} < \frac{3}{4}$.
We find that $h'(y) = \frac{3}{y} - \ln(2 \pi) < 0$ when $y \geq 4$ and that
$\exp(h(4)) = \frac{16}{\pi^3} = \lambda^2$, and it follows that
\begin{eqnarray}
2 y^3 \leq \lambda^2 (2 \pi)^{y-1} && \text{for all } y \geq 4.
\end{eqnarray}

Similarly consider the function $\ol{h} : y \mapsto 5 \ln(y) + \ln(2) - (y-1) \ln(2 \pi)$ for $y > 0$.
We find that $\ol{h}'(y) = \frac{5}{y} - \ln(2 \pi) < 0$ when $y \geq 7$ and that
$\exp(\ol{h}(7)) = \frac{7^5}{32 \pi^6} < 1$, and it follows that
\begin{eqnarray}
2 y^5 \leq (2 \pi)^{y-1} && \text{for all } y \geq 7.
\end{eqnarray}

Then, we set $\bZ(x) = \sqrt{2 \pi} x^{x+1/2} e^{-x}$ for all $x \geq 0$.
We prove below that
\begin{eqnarray}
\frac{(a+b)!}{a! b!} \leq \frac{\bZ(a+b)}{\bZ(a) \bZ(b)} && \text{for all integers } a, b \geq 1.
\end{eqnarray}

\noindent We study the functions $\ds F : (a,b) \mapsto \frac{(a+b)! \bZ(a) \bZ(b)}{\bZ(a+b) a! b!}$
and $\ds G : (a,b) \mapsto \frac{F(a+1,b)}{F(a,b)}$.
We compute that
\begin{align*}
G(a,b) & = \left(\frac{a+b}{a+b+1}\right)^{a+b+1/2}\left(\frac{a+1}{a}\right)^{a+3/2} \geq \ol{G}(a,b)^{2a+b+2} \text{, where} \\
\ol{G}(a,b) & = \frac{2a+b+2}{(a+b+1/2)\frac{a+b+1}{a+b}+(a+3/2)\frac{a}{a+1}} \tag{\text{by geometric-harmonic inequality}} \\
& = 1 + \frac{b-1}{(2a+2b+1)(a+b+1)(a+1)+(2a+3)a(a+b)}.
\end{align*}
and since $b \geq 1$, it follows that $G(a,b) \geq \ol{G}(a,b)^{2a+b+2} \geq 1$, i.e. that $F(a,b) \leq F(a+1,b)$.
Since $F(a,b) = F(b,a)$ for all $a,b \geq 1$, we derive immediately that $F(a,b) \leq F(a,b+1) \leq F(a+1,b+1)$
for all integers $a, b \geq 1$.
Moreover, Stirling's approximation formula states that $a! \sim \bZ(a)$ when $a \to +\infty$.
This proves that $\lim_{\alpha,\beta \to +\infty} F(\alpha,\beta) = 1$, and it follows that $F(a,b) \leq 1$ for all $a, b \geq 1$,
which is indeed equivalent to the inequality (5).

As a corollary, observe that, for all integers $i_1,\ldots,i_m \geq 1$ and using inequality (5), we also have
\begin{align*}
\frac{(i_1+\ldots+i_m)!}{i_1 ! \ldots i_m!} = \prod_{j=2}^m \frac{(i_1+\ldots+i_j)!}{(i_1+\ldots+i_{j-1})!i_j!}
\leq \prod_{j=2}^m \frac{\bZ(i_1+\ldots+i_j)}{\bZ(i_1+\ldots+i_{j-1})\bZ(i_j)},
\end{align*}
from which follows our last auxiliary inequality:
\begin{eqnarray}
\frac{(i_1+\ldots+i_m)!}{i_1 ! \ldots i_m!} \leq \frac{\bZ(i_1+\ldots+i_m)}{\bZ(i_1)\ldots\bZ(i_m)} && \text{for all integers } i_1,\ldots,i_m \geq 1.
\end{eqnarray}

\section{Evaluating $\bK(m,i)$}

We first evaluate $\bT(m,k,\ell)$ when $\ell = 1$.
Here, instead of considering a tuple of non-negative integers $(i_1,\ldots,i_m)$ that sum up to $\ell$,
we might directly consider the unique integer $j \in \{1,\ldots,m\}$ such that $i_j = 1$.
Moreover, for each tuple $(i_1,\ldots,i_m)$, the multinomial coefficient $\binom{\ell}{i_1~\ldots~i_m}$ is equal to $1$.
It follows that $\bT(m,k,1) = m^{-k} \sum_{j=1}^m 1 = m^{1-k}$,
from which we derive the inequalities $\bS(m,k) \geq \bT(m,k,1) = m^{1-k}$ and 
\formula{ \bK(m,i) \leq 2 \prod_{j=1}^{i-1} m^{2^j-1} = 2 \frac{m^{2^i}}{m^{i+1}} }

Then, we investigate lower bounds of $\bK(m,i)$, i.e. upper bounds of $\bT(m,k+1,\ell)$ and of $\bS(m,k+1)$ when $m \geq 4$ and $k \geq 1$.
Consider some integer $\ell \geq 1$, and let us write $\ell = a m + b$, with $a \geq 0$ and $1 \leq b \leq m$.
In addition, let us set
\formula{ \bV(m,a,b) = \max_{i_1+\ldots+i_m=\ell} \binom{\ell}{i_1~\ldots~i_m} \text{ and } \bU(m,a,b) = \frac{1}{m^\ell} \bV(m,a,b)}
We first observe that 
\begin{align*}
\bT(m,k+1,\ell) & = \frac{1}{m^{(k+1)\ell}} \sum_{i_1+\ldots+i_m=\ell} \binom{\ell}{i_1~\ldots~i_m}^{k+1} \\
& \leq \frac{1}{m^{(k+1)\ell}} \sum_{i_1+\ldots+i_m=\ell} \binom{\ell}{i_1~\ldots~i_m} \bV(m,a,b)^k \\
& \leq \bU(m,a,b)^k \bT(m,1,\ell) \\
& \leq \bU(m,a,b)^k. \tag{\text{by Newton multinomial identity}}
\end{align*}

Since the inequality $(x+1)!(y-1)! \geq x! y!$ holds for all integers $x \geq y$,
it also follows that
\formula{ \bU(m,a,b) = \frac{(am+b)!}{m^{am+b} (a+1)!^b a!^{m-b}} }
We compute immediately that $\bU(m,0,b) = m^{-b} b! = m^{-1}$ if $b = 1$, and that $\bU(m,0,b) \leq 2 m^{-2}$ if $2 \leq b \leq m$.
When $a \geq 1$, we further compute that
\begin{align*}
\bU(m,a,b) & \leq \frac{\bZ(am+b)}{m^{am+b} \bZ(a+1)^b \bZ(a)^{m-b}} \tag{\text{using inequality (6)}} \\
& \leq \frac{\sqrt{am+b}}{(2\pi)^{(m-1)/2} a^{m/2}} \frac{(1+b/am)^{am}(1+b/am)^b}{(1+1/a)^{(a+3/2)b}} \\
& \leq \frac{\sqrt{2am}}{(2\pi)^{(m-1)/2} a^{m/2}} \frac{(1+b/am)^{am}(1+1/a)^b}{(1+1/a)^{(a+3/2)b}} & \tag{\text{since $b \leq m \leq am$}} \\
& \leq \frac{\sqrt{2m}}{(2 a \pi)^{(m-1)/2}} \frac{(1+b/am)^{am}}{(1+1/a)^{(a+1/2)b}} \\
& \leq \frac{\sqrt{2m}}{(2 a \pi)^{(m-1)/2}} \frac{e^b}{(1+1/a)^{(a+1/2)b}} & \tag{\text{using inequality (1)}} \\
& \leq \frac{\sqrt{2m}}{(2 a \pi)^{(m-1)/2}} \left(\frac{e}{(1+1/a)^{a+1/2}}\right)^b \\
& \leq \frac{\sqrt{2m}}{(2 a \pi)^{(m-1)/2}}. & \tag{\text{using inequality (2)}}
\end{align*}

Consequently, since $k \geq 1$, we find that
\begin{align*}
\bS(m,k+1) & \leq \sum_{a = 0}^\infty \sum_{b=1}^m \bU(m,a,b)^k
= \bU(m,0,1)^k + \sum_{b=2}^m \bU(m,0,b)^k + \sum_{a=1}^\infty \sum_{b=1}^m \bU(m,a,b)^k \\
& \leq \frac{1}{m^k} + \frac{2(m-1)}{m^{2k}} + m \sum_{a=1}^\infty \left( \frac{\sqrt{2m}}{(2 a \pi)^{(m-1)/2}}\right)^k \\
& \leq \frac{1}{m^k} + \frac{2m}{m^{2k}} + m \frac{(2m)^{k/2}}{(2 \pi)^{k(m-1)/2}} \zeta(k(m-1)/2).
\end{align*}

If we set
\formulaa{ \bP(m,k) = 1 + 2 m^{1-k} + m^{k+1} \frac{(2m)^{k/2}}{(2 \pi)^{k(m-1)/2}} \zeta(k(m-1)/2) }
then it follows that $\bS(m,k+1) \leq \frac{1}{m^k} \bP(m,k)$ and therefore that
\formulaa{ \bK(m,i) \geq 2 \prod_{j=1}^{i-1} \frac{m^{2^j-1}}{\bP(m,2^j-1)} \geq \frac{2}{\bP_\infty(m)} \frac{m^{2^i}}{m^{i+1}} }
where $\bP_\infty(m)$ is the infinite product $\ds \prod_{j=1}^{\infty} \bP(m,2^j-1)$.
It remains to prove that $\bP_\infty(m) \leq 42$.

We fisrt assume that $7 \leq m$. For $k \geq 1$, we compute that
\begin{align*}
\bP(m,k) & = 1 + 2 m^{1-k} + m^{k+1} \frac{(2m)^{k/2}}{(2 \pi)^{k(m-1)/2}} \zeta(k(m-1)/2) \\
& = 1 + 2 m^{1-k} + m^{1-k} \left(\frac{2m^5}{(2\pi)^{m-1}}\right)^{k/2} \zeta(k(m-1)/2) \\
& \leq 1 + 2 m^{1-k} + m^{1-k} \zeta(k(m-1)/2) \tag{\text{using inequality (4)}} \\
& \leq 1 + 4 m^{1-k} \tag{\text{since $\zeta(k(m-1)/2) \leq \zeta(3) \leq 2$}} \\
& \leq \exp(4 m^{1-k}), \tag{\text{since $1+x \leq \exp(x)$ for all $x \in \RR$}}
\end{align*}
from which we deduce that $\bP_\infty(m,1) \leq 5$, whence
\begin{align*}
\ln(\bP_\infty(m)) & = \sum_{j=1}^\infty \ln(\bP(m,2^j-1))
\leq \ln(\bP(m,1)) + 4 \sum_{j=2}^\infty m^{2-2^j}
\leq \ln(5) + 4 \sum_{j=0}^\infty m^{-2-j} \\
& \leq \ln(5) + \frac{4}{m(m-1)}
\leq \ln(5) + \frac{2}{21}
\leq \ln(42). \tag{\text{since $7 \leq m$}}
\end{align*}

Then, we assume that $4 \leq m \leq 6$.
Again, for $k \geq 1$, we compute that
\begin{align*}
\bP(m,k) & = 1 + 2 m^{1-k} + m^{k+1} \frac{(2m)^{k/2}}{(2 \pi)^{k(m-1)/2}} \zeta(k(m-1)/2) \\
& = 1 + 2 m^{1-k} + m \left(\frac{2m^3}{(2\pi)^{m-1}}\right)^{k/2} \zeta(k(m-1)/2) \\
& \leq 1 + 2 m^{1-k} + m \zeta(k(m-1)/2) \lambda^k \tag{\text{using inequality (3)}} \\
& \leq 1 + 2 m^{1-k} + 3 m \lambda^k \tag{\text{$\zeta(k(m-1)/2) \leq \zeta(3/2) \leq 3$}} \\
& \leq 1 + 5 m \lambda^k \tag{\text{since $m^{-1} \leq \frac{1}{4} \leq \lambda$}} \\
& \leq \exp(5 m \lambda^k) \tag{\text{since $1+x \leq \exp(x)$ for all $x \in \RR$}}
\end{align*}
Furthermore, explicit computations in each of the cases $m = 4$, $m = 5$ and $m = 6$ indicate
that $\prod_{j=1}^4\bP(m,2^j-1) \leq 41$. Hence, we conclude that
\begin{align*}
\ln(\bP_\infty(m)) & \leq \ln(41) + \sum_{j=5}^\infty \ln(\bP(m,2^j-1)) \\
& \leq \ln(41) + 5 m \sum_{j=5}^\infty \lambda^{2^j-1}
\leq \ln(41) + 5 m \sum_{j=0}^\infty \lambda^{31+j} \\
& \leq \ln(41) + \frac{5 m \lambda^{31}}{1-\lambda}
\leq \ln(41) + \frac{30 \times 3^{31}}{4^{30}} \leq \ln(42). \tag{\text{since $m \leq 6$ and $\lambda < \frac{3}{4}$}}
\end{align*}

\label{finalpage}


\end{document}